\documentclass[a4paper,12pt]{amsart}
\usepackage{graphicx}
\usepackage{amscd}
\usepackage{amssymb}
\newtheorem{thm}{Theorem}[section]
\newtheorem{cor}[thm]{Corollary}
\newtheorem{lem}[thm]{Lemma}
\newtheorem{prop}[thm]{Proposition}
\theoremstyle{defn}
\newtheorem{defn}[thm]{Definition}
\theoremstyle{rem}
\newtheorem{rem}[thm]{Remark}

\numberwithin{equation}{section}

\newcommand{\eps}{\varepsilon}

\newcommand{\G}{\mathcal{G}}

\newcommand{\T}{\mathcal{T}(\mathcal G)}

\newcommand{\Euv}{\mathcal E_{u,v}^\pi}
\newcommand{\R}{\mathcal{R}ep(\mathcal G)}
\newcommand{\F}{\mathfrak {F}_{u,v}}
\newcommand{\HH}{\mathcal H_\pi}
\newcommand{\Hu}{\mathcal H_u^\pi}
\newcommand{\Hv}{\mathcal H_v^\pi}

\newcommand{\repn}{representation\,}
\newcommand{\repns}{representations\,}

\begin{document}

\title[duality for groupoids]{Tannaka-Krein duality for compact groupoids II,
Fourier transform}%
\author{Massoud Amini}%
\address{
Department of Mathematics, Tarbiat Modarres University, P.O.Box
14115-175, Tehran , Iran , mamini@modares.ac.ir \newline
Department of Mathematics and Statistics, University of
Saskatchewan, 106 Wiggins Road, Saskatoon, S7N 5E6 ,
mamini@math.usask.ca}

\thanks{The author was visiting the University of Saskatchewan during the
preparation of this work. He would like to thank University of
Saskatchewan and Professor Mahmood Khoshkam for their hospitality
and support}%


\begin{abstract}
In a series of papers, we have shown that from the \repn theory of
a compact groupoid one can reconstruct the groupoid using the
procedure similar to the Tannaka-Krein duality for compact groups.
In this part we study the Fourier and Fourier-Plancherel
transforms and prove the Plancherel theorem for compact groupoids.
We also study the central functions in the algebra of square
integrable functions on the isotropy groups.

\end{abstract}
\maketitle
\section{introduction}
In a series of papers, we have generalized the Tannaka-Krein
duality to compact groupoids. In [A1] we studied the \repn theory
of compact groupoids. In particular, we showed that irreducible
\repns have finite dimensional fibres. We also proved the Schur's
lemma, Gelfand-Raikov theorem and Peter-Weyl theorem for compact
groupoids. In this part we study the Fourier and
Fourier-Plancherel transforms on compact groupoids. In section two
we develop the theory of Fourier transforms on the Banach algebra
bundle $L^1(\G)$ of a compact groupoid $\G$. As in the group case,
a parallel theory of Fourier-Plancherel transform on the Hilbert
space bundle $L^2(\G)$ is constructed. This provides a surjective
isometric linear isomorphism from $L^2(\G)$ to $L^2(\hat \G)$, in
an appropriate sense. Also the relation between $\hat \G$ and the
conjugacy groupoid $\G^\G$ is studied. The results of this section
are effectively used in [A2] to show that the natural homomorphism
from $\G$ to its Tannaka groupoid $\T$ is surjective. Section
three considers the inverse Fourier and Fourier-Plancherel
transforms. In this section we prove Plancherel theorem for
compact groupoids. In section four we study the class functions
and central elements in the algebras of functions on fibres of
$\G$ and prove a diagonal version of the Plancherel theorem. All
over this paper we assume that $\G$ is compact and the Haar system
on $\G$ is normalized.

\vspace{.5 cm}
\section{Fourier transform}

 It follows from the Peter-Weyl theorem [A1, theorem 3.13] that, for $u,v\in X$,
 if $f\in L^2(\G_u^v,\lambda_u^v)$ then
$$f=\sum_{\pi\in\hat\G}\sum_{i=1}^{d_v^\pi}\sum_{j=1}^{d_u^\pi}
c_{u,v,\pi}^{ij} \pi_{u,v}^{ij}\,\,,$$ where
$$c_{u,v,\pi}^{ij}=d_u^\pi\int_{\G_u^v} f(x)
\overline{\pi_{u,v}^{ij}(x)}d\lambda_u^v(x)\quad(1\leq i\leq
d_v^\pi, \,1\leq j\leq d_u^\pi).$$

This is a local version  of the classical non commutative Fourier
transform. As in the classical case, the main drawback is that it
depends on the choice of the basis (which in turn gives the choice
of the coefficient functions). The trick is similar to the
classical case, that's to use the continuous decomposition using
integrals. This is the content of the next definition. As usual,
all the integrals are supposed to be on the support of the measure
against which they are taken.

\begin{defn}
Let $u,v\in X$ and $f\in L^1(\G_u^v,\lambda_u^v)$, then the
Fourier transform of $f$ is $\F(f):\R\to \mathcal B(\Hv,\Hu)$
defined by
$$\F(f)(\pi)=\int f(x)\pi(x^{-1})d\lambda_u^v(x).$$
\end{defn}

To  better understand this definition, let us go back to the group
case for a moment. Let's start with a locally compact Abelian
group $G$. Then the Pontryagin dual $\hat G$ of $G$ is a locally
compact Abelian group and for each $f\in L^1(G)$, its Fourier
transform $\hat f\in C_0(\hat G)$ is defined by
$$\hat f(\chi)=\int_G f(x)\overline{\chi(x)} dx\quad (\chi\in \hat G).$$
The continuity of $\hat f$  is immediate and the fact that it
vanishes at infinity is the so called {\it Riemann-Lebesgue
lemma}. In the non Abelian compact case, one get a similar
construction, namely, with an slight abuse of notation, for each
$f\in L^1(G)$ one has $\hat f\in C_0(\hat G,\mathcal B(\mathcal
H))$, where $\hat G$ is the set of (unitary equivalence classes
of) irreducible \repns of $G$ endowed with the Fell topology, and
$\mathcal B(\mathcal H)$ is a bundle of $C^{*}$-algebras over
$\hat G$ whose fiber at $\pi$ is $\mathcal B(\mathcal H_\pi)$, and
by $C_0(\hat G,\mathcal B(\mathcal H))$ we mean the set of all
continuous sections which vanish at infinity. In the groupoid
case, one has a similar interpretation. Locally each $f\in
L^1(\G_u^v,\lambda_u^v)$ has its Fourier transform $\F(f)$ in
$C_0(\hat \G,\mathcal B_{u,v}(\mathcal H))$, where $\hat \G$ is
the set of (unitary equivalence classes of) irreducible \repns of
$\G$ endowed again with the Fell topology, and $\mathcal
B_{u,v}(\mathcal H)$ is a bundle of $C^{*}$-algebras over $\hat
\G$ whose fiber at $\pi$ is $\mathcal B(\Hv,\Hu)$, and $C_0(\hat
G,\mathcal B_{u,v}(\mathcal H))$ is the set of all continuous
sections vanishing at infinity. Globally we have an still more
complicated interpretation. We have to look at $L^1(\G)$ as a
bundle of Banach algebras over $\G^{(0)}\times\G^{(0)}$ whose
fiber at $(u,v)$ is $L^1(\G_u^v,\lambda_u^v)$, and then each $f\in
L^1(\G)$ has its Fourier transform $\mathfrak F(f)$ in $C_0(\hat
\G,\mathcal B(\mathcal H))$, where $\mathcal B(\mathcal H)$ is a
bundle of bundles of $C^{*}$-algebras over $\hat \G$ whose fiber
at $\pi$ is the bundle $\mathcal B(\HH)$ over
$\G^{(0)}\times\G^{(0)}$ whose fiber at $(u,v)$ is $\mathcal
B(\Hv,\Hu)$, the space $C_0(\hat G,\mathcal B(\mathcal H))$ is the
set of all continuous sections vanishing at infinity, and
$\mathfrak {F}(f)(\pi)_{(u,v)}=\F(f_{(u,v)})(\pi)$!!

Now  let us discuss the properties of the Fourier transform. If we
choose (possibly infinite) orthonormal bases for $\Hu$ and $\Hv$
and let each $\pi(x)$ be represented by the (possibly infinite)
matrix with components $\pi_{u,v}^{ij}(x)$, then $\F(f)$ is
represented by the matrix with components
$\F(f)(\pi)^{ij}=\frac{1}{d_u^\pi}c_{u,v,\pi}^{ji}$\,. When $f\in
L^2(\G_u^v,\lambda_u^v)$, summing up over all indices $i,j$, we
get the following

\begin{prop}{\bf (Fourier inversion formula)}
For each $u,v\in X$ and $f\in L^2(\G_u^v,\lambda_u^v)$,
$$f=\sum_{\pi\in\hat\G} d_u^\pi Tr\big(\F(f)(\pi)\pi(.)\big),$$
where the sum converges in the $L^2$ norm and
$$\|f\|_2^2=\sum_{\pi\in\hat\G} d_u^\pi Tr\big(\F(f)(\pi)^{*}\F(f)(\pi)\big) .$$\qed
\end{prop}

We collect  the properties of the Fourier transform in the
following lemma. The proof is routine and is omitted.

\begin{lem}
Let $u,v\in X$, $a,b\in\mathbb C$, $x\in\G$, and $f,g\in
L^1(\G_u^v,\lambda_u^v)$, then for each $\pi\in\R$,

$(i) \F(af+bg)=a\F(f)+b\F(g)$,

$(ii) \F(f*g)(\pi)=\F(f)(\pi)\F(g)(\pi)$,

$(iii) \F(f^*)(\pi)=\F(f)(\pi)^*$,

$(iv) \F(\ell_x(f))(\pi)=\F(f)(\pi)\pi(x^{-1}), \, \,
\F(r_x(f))(\pi)=\pi(x)\F(f)(\pi)\quad (x\in\G_u^v)$.\qed
\end{lem}

\begin{cor}
 For $a,b\in\mathbb C$, $f,g\in L^1(\G)$, and $\pi\in\hat \G$,

$(i) \mathfrak F(af+bg)=a\mathfrak F(f)+b\mathfrak F(g)$,

$(ii) \mathfrak F(f*g)(\pi)=\mathfrak F(f)(\pi)\mathfrak
F(g)(\pi)$,

$(iii) \mathfrak F(f^*)(\pi)=\mathfrak F(f)(\pi)^*$,

\end{cor}

As in  the group case there is yet another way of introducing the
Fourier transform. For each finite dimensional continuous \repn
$\pi$ of $\G$, let the {\it character} $\chi_\pi$ of $\pi$ be the
bundle of functions $\chi_\pi$ whose fiber at $u\in X$ is
$\chi_u^\pi(x)=Tr(\pi(x)) \quad (x\in \G_u^u)$, where $Tr$ is the
trace of matrices. Note that one can not have these as functions
defined on $\G_u^v$, since when $x\in\G_u^v$,  $\pi(x)$ is not a
square matrix in general. Also note that the values of the above
character functions depend only on the unitary equivalence class
of $\pi$, as similar matrices have the same trace. Now if
$\pi\in\hat\G$,  $x\in\G_u^v$, and $f\in L^1(\G_u^v,\lambda_u^v)$,
then
$$Tr\big(\F(f)(\pi)\pi(x)\big)=\int f(y)Tr(\pi(y^{-1}x))d\lambda_u^v(y)=f*\chi_u^\pi(x),$$
so it follows from Proposition 2.2 that

\begin{cor}
The map $P_{u,v}^\pi: L^2(\G_u^v,\lambda_u^v)\to \Euv, \quad
f\mapsto d_u^\pi f*\chi_u^\pi$ is a  surjective orthogonal
projection and for each $f\in L^2(\G_u^v,\lambda_u^v)$ we have the
decomposition
$$f=\sum_{\pi\in\hat\G} d_u^\pi f*\chi_u^\pi,$$
which converges in the $L^2$ norm.\qed
\end{cor}

Applying the above decomposition to the case where $u=v$ and
$f=\chi_u^\pi$, we get

\begin{cor}
For each $u\in X$ and $\pi,\pi^{'}\in\hat \G$,
\begin{equation*}
\chi_u^\pi*\chi_u^{\pi^{'}}=
\begin{cases}
\, {d_u^\pi}^{-1}& \text{if} \,\, \pi\sim\pi^{'},\\
\,0& \text{otherwise}.
\end{cases}
\end{equation*}
\qed
\end{cor}

\vspace{.5 cm}
\section{inverse Fourier and Fourier-Palncherel transforms}

Next we are aiming at the construction  of the inverse Fourier
transform. This is best understood if we start with a yet
different interpretation of the local Fourier transform. It is
clear from the definition that if $u,v\in X$, $\pi_1,\pi_2\in\R$,
and $f\in L^1(\G_u^v,\lambda_u^v)$, then
$$\F(f)(\pi_1\oplus\pi_2)=\F(f)(\pi_1)\oplus\F(f)(\pi_2),$$
and the same is true for any number  (even infinite) of continuous
\repns , so it follows from Theorem 2.16 in [A1] that $\F(f)$ is
uniquely characterized by its values on $\hat\G$, namely we can
regard
$$\F :L^1(\G_u^v,\lambda_u^v)\to\prod_{\pi\in\hat\G} \mathcal B(\Hv,\Hu).$$
Now consider the $C^*$-algebra $\ell^\infty$-direct  sum
$\sum_{\pi\in\hat\G} \bigoplus \mathcal B(\Hv,\Hu)$. The domain of
our inverse Fourier transform then would be the algebraic sum
$\ell^\infty$-direct sum $\sum_{\pi\in\hat\G} \mathcal
B(\Hv,\Hu)$, consisting of those elements of this $C^*$-algebra
with only finitely many nonzero components.

\begin{defn}
Let $u,v\in X$. The inverse Fourier transform
$$\F^{-1}:\sum_{\pi\in\hat\G} \mathcal B(\Hv,\Hu)\to C(\G_u^v)$$
is defined by
$$\F^{-1}(g)(x)= \sum_{\pi\in\hat\G} d_u^\pi Tr\big(g(\pi)\pi(x^{-1})\big)
\quad (x\in\G_u^v).$$
\end{defn}

To show that this is indeed the inverse map of  the (local)
Fourier transform we need some orthogonality relations. They are a
version of the Schur's orthogonality relations [A1, theorem 3.6]
and proved similarly, so We only give a sketch of the proof.

\begin{prop} {\bf (Orthogonality relations)}
Let $\tau,\rho\in\hat \G$, $u,v\in X$, $T\in\mathcal B(\mathcal
H_\tau)$, $S\in \mathcal B(\mathcal H_\rho)$, $\xi\in\mathcal
H_\tau$, $\eta\in\mathcal H_\rho$, and $A\in Mor(\tau,\rho)$, then

\begin{equation*}
 (i)  \quad\int \tau(x)A_{s(x)}\rho(x^{-1})d\lambda_u(x)=
\begin{cases}
\, \frac{Tr(A_u)}{d_u^\tau} id_{\mathcal H_u^\tau}& \text{if} \,\, \tau=\rho,\\
\,0& \text{otherwise},
\end{cases}
\end{equation*}

\begin{equation*}
(ii)\quad\int
\tau(x)\xi_{s(x)}\otimes\rho(x^{-1})\eta_{r(x)}d\lambda_u(x)=
\begin{cases}
\, \frac{\eta_u\otimes\xi_u}{d_u^\tau} & \text{if} \,\, \tau=\rho,\\
\,0& \text{otherwise},
\end{cases}
\end{equation*}

\begin{equation*}
(iii)\quad\int
Tr(T_{r(x)}\tau(x))Tr(S_{s(x)}\rho(x^{-1})d\lambda_u(x)=
\begin{cases}
\, \frac{Tr(T_uS_u)}{d_u^\tau} & \text{if} \,\, \tau=\rho,\\
\,0& \text{otherwise},
\end{cases}
\end{equation*}

\begin{equation*}
(iv)\quad\int
Tr(T_{r(x)}\tau(x))\overline{Tr(S_{r(x)}\rho(x)}d\lambda_u(x)=
\begin{cases}
\, \frac{Tr(T_uS_u^{*})}{d_u^\tau} & \text{if} \,\, \tau=\rho,\\
\,0& \text{otherwise},
\end{cases}
\end{equation*}

\begin{equation*}
(v)\quad\int Tr(T_{r(x)}\tau(x))\rho(x^{-1})d\lambda_u(x)=
\begin{cases}
\, \frac{1}{d_u^\tau}T_u & \text{if} \,\, \tau=\rho,\\
\,0& \text{otherwise}.
\end{cases}
\end{equation*}

\end{prop}
{\bf Proof} $(i)$ As in  [A1, lemma 3.4], the left hand side
defines a bundle of operators in $Mor(\tau,\rho)$, so by Schur's
lemma [A1, theorem 2.14] it is $c.id_{\mathcal H_u^\tau}$, if
$\tau=\rho$, and $0$, otherwise. Now
$$Tr\big(\int \tau(x)A_{s(x)}\tau(x^{-1})d\lambda_u(x)\big)=Tr(A_u),$$
where as $Tr\big(c.id_{\mathcal H_u^\tau}\big)=cd_u^\tau$, so $c$
is what it should be.

$(ii)$ Take any $\phi,\psi\in\mathcal H_\tau^{*}$  and apply $(i)$
to $A$ defined by $A_u(\zeta_u)=\phi_u(\zeta_u)\xi_u \quad (u\in
X)$ and then calculate both sides of the resulting operator
equation at $\eta_u$ to get

\begin{equation*}
\int \tau(x)\xi_{s(x)}\phi(\rho(x^{-1}\eta_u)d\lambda^u(x)=
\begin{cases}
\, \frac{\phi_u(\xi_u)}{d_u^\tau}\eta_u & \text{if} \,\, \tau=\rho,\\
\,0& \text{otherwise}.
\end{cases}
\end{equation*}

The result now follows if we apply $\psi_u$ to both  sides of the
above equality and use the fact that
$\phi_u(\xi_u)\psi_u(\eta_u)=(\psi\otimes\phi)_u(\eta_u\otimes\xi_u)$.

$(iii),(iv)$ Note that all the involved Hilbert spaces  are finite
dimensional [A1, theorem 2.16]. In particular, rank one operators
generate all operators on these spaces. Also the required relation
is linear in $T$ and $S$. Hence we may assume that $T$ and $S$
have rank one fibers, say $T=\phi(.)\xi, S=\psi(.)\eta$, where
$\phi,\psi$ are as above. Now applying $(\phi\otimes\psi)_u$ to
both sides of $(ii)$, we get $(iii)$. The proof of $(iv)$ is
similar.

$(v)$ Let $L$ and $R$ be the left and right hand sides of  $(v)$,
respectively. We need only to show that $Tr((L-R)S)=0$, for each
$S\in\mathcal B(\mathcal H_\rho)$. But $Tr(LS)$ is clearly the
right hand side of $(iii)$, which is in turn equal to
$Tr(RS)$.\qed

\vspace{.3 cm} Now we are ready to prove the properties of the
local inverse  Fourier transform. But let us first introduce the
natural inner products on its domain and range. For $f,g\in
C(\G_u^v)$ and $h,k\in\sum_{\pi\in\hat\G} \mathcal B(\Hv,\Hu)$ put
$$<f,g>=\int \bar f .g d\lambda_u^v,$$
and
$$<h,k>=\sum_{\pi\in\hat\G} d_u^\pi Tr(h^*(\pi)k(\pi)),$$
where the right hand  side is a finite sum as $h$ and $k$ are of
finite support. Also note that if $\eps_u:C(\G_u^u)\to \mathbb C$
is defined by $\eps_u(f)=f(u)$, then for each $f,g\in C(\G_u^v)$,
we have $f^**g\in C(\G_u^u)$ and $<f,g>=\eps_u(f^**g)$.

\begin{prop}
For each $u,v\in X$ and $h,k\in\sum_{\pi\in\hat\G} \mathcal
B(\Hv,\Hu)$ we have

$(i) \F\F^{-1}(h)=h$,

$(ii) \F^{-1}(hk)=\F^{-1}(h)*\F^{-1}(k)$,

$(iii) \F^{-1}(h^*)=(\F^{-1}(h))^*$,

$(iv) <\F^{-1}(h),\F^{-1}(k)>=<h,k>$.
\end{prop}
{\bf Proof} $(i)$ By $(v)$ of above proposition, for each
$\tau\in\hat\G$,
$$\F\F^{-1}(h)(\tau)=\sum_{\pi\in\hat\G} \int d_u^\pi Tr(h(\pi)\pi(x^{-1}))
\tau(x)d\lambda_u^v(x)=h(\tau).$$

$(ii)$ By $(iii)$ of above proposition, for each $x\in\G_u^v$,
\begin{align*}
(\F^{-1}(h)*\F^{-1}(k))(x)&=\int \F^{-1}(h)(xy^{-1})\F^{-1}(k)(y) d\lambda_u^v(y)\\
&=\sum_{\tau,\rho\in\hat\G} \int d_u^\tau d_u^\rho
Tr\big(h(\tau)\tau(yx^{-1})\big) Tr\big(k(\rho)\rho(y^{-1})\big)d\lambda_u^v(x)\\
&=\sum_{\tau,\rho\in\hat\G} d_u^\tau d_u^\rho\int
Tr\big(h(\tau)\tau(x^{-1})\tau(y)\big)
 Tr\big(k(\rho)\rho(y^{-1})\big)d\lambda_u^v(x)\\
&=\sum_{\tau\in\hat\G} d_u^\tau   Tr\big(h(\tau)\tau(x^{-1})k(\tau)\big)\\
&=\sum_{\tau\in\hat\G} d_u^\tau   Tr\big(h(\tau)k(\tau)\tau(x^{-1})\big)\\
&=\F^{-1}(hk)(x).
\end{align*}

$(iii)$ For each $x\in\G_u^v$
\begin{align*}
\F^{-1}(h^{*})(x)&=\sum_{\tau\in\hat\G} d_u^\tau   Tr\big(h^{*}(\tau)\tau(x^{-1})\big)\\
&=\sum_{\tau\in\hat\G} d_u^\tau   Tr\big(\bar h((\check\tau)^{\bar{}})\tau(x^{-1})\big)\\
&=\sum_{\tau\in\hat\G} d_u^\tau   Tr\big(\bar h(\tau)\overline{\check\tau}(x^{-1})\big)\\
&=\sum_{\tau\in\hat\G} d_u^\tau  Tr\big(\bar
h(\tau)\bar\tau(x)\big),
\end{align*}
where as
\begin{align*}
(\F^{-1}(h))^{*}(x)&=\overline{\F^{-1}(h)(x^{-1})}\\
&=\sum_{\tau\in\hat\G} d_u^\tau   \overline{Tr\big(h(\tau)\tau(x)\big)}\\
&=\sum_{\tau\in\hat\G} d_u^\tau   Tr\big(\bar
h(\tau)\bar\tau(x)\big).
\end{align*}

$(iv)$ By above observation about $\eps_u$,
\begin{align*}
<\F^{-1}(h),\F^{-1}(k)>&=\eps_u((\F^{-1})^**\F^{-1}(h))=\eps_u(\F^{-1}(h^*k))
=\F^{-1}(h^*k)(u)\\
&=\sum_{\tau\in\hat\G} d_u^\tau
Tr\big(h^*(\tau)k(\tau)\big)=<h,k>.
\end{align*}
\qed

\vspace{.3 cm} Next we define a norm on the domain of the inverse
Fourier transform in order to get a Plancherel type theorem. For
$h\in \sum_{\pi\in\hat\G} \mathcal B(\Hv,\Hu)$ we put
$\|h\|_2=<h,h>^{\frac{1}{2}}$. This is the natural norm on the
algebraic direct sum, when one endows each component $\mathcal
B(\Hv,\Hu)$  with the pre-Hilbert space structure given by
$<T,S>=\big( d_u^\pi Tr(S^*T)\big)^{\frac{1}{2}}$. We denote the
completion of $\sum_{\pi\in\hat\G} \mathcal B(\Hv,\Hu)$ with
respect to this norm by $\mathcal L_{u,v}^2(\G)$. We deliberately
used the curly $\mathcal L$ for this space to distinguish it from
its counterpart which is defined later.

\begin{thm} {\bf (Plancherel Theorem)}
For each $u,v\in X$, $\F$ extends (uniquely) to a continuous
surjective linear isometry
 $\F: L^2(\G_u^v,\lambda_u^v)\to\mathcal L_{u,v}^2(\G)$.
\end{thm}
{\bf Proof} By Proposition 3.4 $(iv)$, $\F^{-1}:\mathcal
L_{u,v}^2(\G)\to L^2(\G_u^v,\lambda_u^v)$  is an isometric
embedding. It is also surjective, since $Im(\F^{-1})$ is complete
and so closed, and also it clearly includes $\mathcal E_{u,v}$
which is dense in $L^2(\G_u^v,\lambda_u^v)$. \qed

\vspace{.3 cm} The above map is called the (local){\it
Fourier-Plancherel transform}. Now for each $u,v\in X$, $\G_u^u$
and $\G_v^v$  act on both $L^2(\G_u^v,\lambda_u^v)$ and $\mathcal
B(\Hv,\Hu)$, from right and left respectively, via
$$(f.x)(y)=f(yx),\, A.x=\pi(x)A\quad (x\in\G_u^u, y\in\G_u^v),$$
and
$$(x.f)(y)=f(x^{-1}y),\, x.A=A\pi(x)\quad (x\in\G_v^v, y\in\G_u^v),$$
for each $f\in L^2(\G_u^v,\lambda_u^v), A\in \mathcal B(\Hv,\Hu)$.

It is easy to see that the Fourier-Plancherel transform respects
these actions, namely

\begin{lem}
For each $u,v\in X$, $x\in\G_v^v$, $y\in\G_u^u$, and $f\in
L^2(\G_u^v,\lambda_u^v)$,
$$\F(x.f)=x.\F(f),\, \F(f.y)=x.\F(f).$$\qed
\end{lem}

Before we end this section, let us show that  how one can use
characters of \repns in $\hat \G$ and the orthogonality relations
of the beginning of this section to prove statements about subsets
of $\R$.

\begin{lem}
For each $u\in X$ and $\pi\in\hat\G$, $\chi_u^\pi\in\mathcal
E_{u,u}$.
\end{lem}
{\bf Proof} Let $\{e_u^i\}_{1\leq i\leq d_u^\pi}$ be a basis for
$\Hu$, then
$$\chi_u^\pi=Tr(\pi(.))=\sum_{i=1}^{d_u^\pi} <\pi(.)e_u^\pi,e_u^\pi>,$$
so the result follows from Proposition 3.2 of [A1].\qed

\begin{defn}
A subset $\Sigma$ of $\R$ is called closed if it contains

$(i)\, \pi_1$ if $\pi_1$ is unitary equivalent to some
$\pi_2\in\Sigma$,

$(ii)\, \pi_1$ if $\pi_1$ is weakly contained in some
$\pi_2\in\Sigma$,

$(iii)\, \pi_1\oplus\pi_2$ if $\pi_1,\pi_2$ are in $\Sigma$,

$(iv)\, \pi_1\otimes\pi_2$ if $\pi_1,\pi_2$ are in $\Sigma$,

$(v)\, \bar\pi_1$ if $\pi_1$ is in $\Sigma$,

$(vi)$ the trivial \repn \, $tr$.
\end{defn}

\begin{prop}
If $\Sigma\subseteq\R$ is closed and separates the points of $\G$,
then $\Sigma=\R$.
\end{prop}
{\bf Proof} If not, by  condition $(iii)$ of the definition of
closedness and Theorem 2.16 of [A1], there is $\tau\in\hat\G$
which is not in $\Sigma$. Let $\mathcal
E_{u,v}^\Sigma=\cup_{\pi\in\Sigma} \Euv$, for $u,v\in X$. By
Proposition 3.2, elements of each $\mathcal E_{u,v}^\tau$ is
orthogonal to $\mathcal E_{u,v}^\Sigma$. In particular, by above
lemma, $\chi_u^\tau\in(\mathcal E_{u,v}^\Sigma)^\perp$. But by
conditions $(iii)$-$(v)$ of the definition of closedness,
$\mathcal E_{u,v}^\Sigma$ is a subalgebra of $C(\G)$ which is
closed under conjugation, and by condition $(vi)$, it contains the
constants, and finally by assumption,it separates the points of
$\G$. Hence, by Stone-Weierstrass Theorem, $\mathcal
E_{u,v}^\Sigma$ is dense in $C(\G)$. Therefore $\chi_u^\tau$ is
orthogonal to $C(\G)$ and so it is zero, which is a
contradiction.\qed

\section{Class functions and central elements}

We characterize central elements in $C(\G_u^v)$ and
$L^2(\G_u^v,\lambda_u^v)$, with respect  to the convolution
product.

\begin{defn}
A function $f\in C(\G)$ is called a class  function if it is
constant on "conjugacy classes" of $\G$, that is
$$f(x^{-1}yx)=f(y)\quad(u\in X, x\in\G^u, y\in\G_u^u).$$
We denote the set of all class functions on $\G$ by $\mathfrak
CC(\G)$.
\end{defn}

When we talk about $C(\G)$ as an algebra we always consider it
with  the pointwise multiplication. However, it is clear that
$C(\G)$ is also an algebra with respect the convolution. This
later algebra is non commutative in general, and to distinguish it
from $C(\G)$ we will denote it with $(C(\G),*)$. Also we denote
the center of and algebra $\mathcal A$ with $\mathcal Z(\mathcal
A)$.

\begin{lem}
$\mathfrak CC(\G)=\mathcal Z(C(\G),*)$.
\end{lem}
{\bf Proof} If $f\in\mathfrak CC(\G)$ and $g\in C(\G)$, then for
each $x\in\G$,
\begin{align*}
(f*g)(x)&=\int f(xy^{-1})g(y)d\lambda_{s(x)}(y)\\
&=\int f(xy^{-1}x^{-1})g(xy)d\lambda^{s(x)}(y)\\
&= \int f(y^{-1})g(xy)d\lambda^{s(x)}(y)\\
&=\int f(y)g(xy^{-1})d\lambda_{s(x)}(y).\\
&=(g*f)(x)
\end{align*}
Conversely  if $f\in \mathcal Z(C(\G),*)$, then for each $u\in X$,
$x\in \G^u$, and $g\in C(\G)$,
$$\int (f(xyx^{-1})-f(y))g(xy^{-1})d\lambda_u^u(y)=0,$$
so by continuity of $f$, $f(xyx^{-1})-f(y)=0$ for each
$y\in\G_u^u$.\qed

\vspace{.3 cm} Recall that  each $f\in L^1(\G)$ has a global
Fourier transform $\mathfrak F(f)$ such that for each
$\pi\in\hat\G$, $\mathfrak F(f)(\pi)$ is fibred over $X\times X$.
For some technical reasons, sometimes we need to consider
everything to be fibred over $X$ (rather than $X\times X$. This is
naturally done by considering a diagonal version of this global
Fourier transform, namely $\mathfrak F(f)(\pi)$ is fibred over $X$
and its fiber at $u$ is $\mathfrak F_{u,u}(f_{(u,u)})$. To
distinguish these two, we denote the diagonal version of the
global Fourier transform by $\mathfrak {DF}$.

Now each $f\in C(\G)$ could be considered as a an element  in
$L^1(\G)$ whose fiber at $(u,v)$ is the restriction of $f$ to
$\G_u^v$.

\begin{lem}
If $f\in\mathfrak CC(\G)$ then for each $\pi\in\hat\G$, $\mathfrak
{DF}(f)(\pi)\in Mor(\pi,\pi)$.
\end{lem}
{\bf Proof} Let $\pi\in\hat\G$. By Lemma 4.2 and  definition of
$\mathfrak CC(\G)$, for each $u\in X$ and $x\in\G_u^u$ we have
$$x.\mathfrak F_{u,u}(f)(\pi)=\mathfrak F_{u,u}(x.f)(\pi)
=\mathfrak F_{u,u}(f.x)(\pi)=\mathfrak F_{u,u}(f)(\pi).x.$$\qed

\vspace{.3 cm} Up  to now our Fourier transform was an operator
valued map. As in the group case, we can use the characters of
irreducible \repns to define a complex valued version of the
Fourier transform.

\begin{defn}
Let $\G$ be a compact  groupoid and $f\in L^1(\G)$. Consider
$L^1(\G)$ as a bundle over $X$ whose fiber at $u\in X$ is
$L^1(\G_u,\lambda_u)$. Define $\hat f$ on $\hat \G$ fibrewise by
$$\hat f_u(\pi)=\frac{1}{d_u^\pi}\int f_u(x)\chi_u^\pi(x^{-1}) d\lambda_u^u(x)
=\frac{1}{d_u^\pi}\int f_u(x)Tr(\pi(x^{-1}))d\lambda_u^u(x)\quad
(u\in X).$$ We  call $\hat f$ the diagonal Fourier transform of
$f$ (the terminology is justified with the next proposition).
\end{defn}

Because of the  restriction imposed by the trace,  we had to look
at the $L^1(\G)$ as a bundle over $X$ (rather than $X\times X$).
As one might guess, this makes our new version of the Fourier
transform compatible with the diagonal form of the global Fourier
transform.

\begin{prop}
For $f\in\mathfrak CC(\G)$, \,\,$\mathfrak {DF}(f)(\pi) =\hat
f(\pi)id_{\HH}\quad (\pi\in\hat\G)$.
\end{prop}
{\bf Proof} First note that the above equality means that for each
$u\in X$, $\mathfrak F_{u,u}(f)(\pi)=\hat f_u(\pi)id_{\Hu}$. It
follows from above lemma and Schur's lemma [A1, theorem 2.14],
that $\mathfrak {DF}(f)(\pi)=c_\pi id_{\HH}$, for some bundle of
constants $c_\pi=\{c_u^\pi\}$. But then for $u\in X$,
$$c_u^\pi=\frac{1}{d_u^\pi} Tr\big(\mathfrak F_{u,u}(f)(\pi)\big)
=\frac{1}{d_u^\pi} Tr\big(\int
f_u(x)\pi(x^{-1})d\lambda_u^u(x)\big)=\hat f_u(\pi).$$\qed

\vspace{.3 cm} The diagonal Fourier transform has the same
properties as the global Fourier transform. The proof is routine.

\begin{prop}
 For $a,b\in\mathbb C$, $f,g\in L^1(\G)$, and $\pi\in\hat \G$,

$(i)(af+bg)\hat{}=a\hat f+b\hat g$,

$(ii) (f*g)\hat{}(\pi)=\hat f(\pi)\hat g(\pi)$,

$(iii) (f^*)\hat{}(\pi)=\hat f(\pi)^*$.\qed

\end{prop}

\begin{rem} In above  proposition, if $f,g\in\mathfrak CC(\G)$,
one can easily check that $af+bg, f*g, f^*\in\mathfrak CC(\G)$. In
this case the above relations follow from Lemmas 2.3 and 4.2.
\end{rem}

Next we turn into the inverse of the diagonal Fourier transform.
Here $\hat \G$ is endowed with the discrete topology, so for
instance $C_c(\hat \G)$ simply means all complex valued functions
on $\hat \G$ with finite support [consider it as a bundle with
fiber at u to be ?]. Also we consider $C(\G)$ as a bundle over $X$
whose fiber at $u\in X$ is $C(\G_u^u)$.

\begin{defn}
Let $\G$ be a compact groupoid. The inverse diagonal Fourier
transform  from $C_c(\hat\G)$ to $C(\G)$ is defined by $g\mapsto
\check g$, where
\begin{equation*}
\check g(x)=
\begin{cases}
\, \sum_{\pi\in\hat\G} d_u^\pi g(\pi)\overline{\chi_u^\pi(x)}&
\text{if} \,\,x\in \G_u^u \,\,\text{for some}\,\, u\in X,\\
\,0& \text{otherwise}.
\end{cases}
\end{equation*}
\end{defn}

We  have already used the notation $\check g$ with a different
meaning, namely for a function $g$ on $\G$, $\check
g(x)=g(x^{-1})$. However, since this notation is now used only for
functions on $\hat \G$, there is no fear of confusion.

\begin{lem}
If $g\in C_c(\hat \G)$  then $\check g\in \mathfrak C(\G)$. In
particular, if $\check g$ is continuous then $\check g\in\mathfrak
CC(\G)$. This is the case, for instance when $X=\G^{(0)}$ is
discrete in the relative topology of $\G$ (when $\G$ is Hausdorff,
this means that $X$ is finite).
\end{lem}
{\bf Proof} First let us prove the last  statement. If
$x_\alpha\to x$ in $\G$, then by continuity of the source and
range maps, $s(x_\alpha)\to u=s(x)$ and $r(x_\alpha)\to v=r(x)$.
If $X$ is discrete, eventually $s(x_\alpha)=u$ and
$r(x_\alpha)=v$. If $u\neq v$, eventually $\check
g(x_\alpha)=\check g(x)=0$, otherwise eventually $\check
g(x_\alpha)=\check g_u(x_\alpha)\to \check g_u(x)=\check g(x)$, by
the fact that the character $\chi_u^\pi$ of each $\pi\in\hat\G$ is
continuous.

If $g\in C_c(\hat\G)$, then for each $f\in C(\G)$ and $x\in \G$
with $u=s(x), v=r(x)$ we have
\begin{align*}
(\check g*f)(x)&=\int \check g(xy^{-1})f(y)d\lambda_{s(x)}(y)\\
&=\int \check g(xy^{-1})f(y)d\lambda_{u}^v(y)\\
&=\sum_{\pi\in\hat\G} \int d_u^\pi g(\pi)\overline{\chi_v^\pi(xy^{-1})} f(y)d\lambda_{u}^v(y)\\
&=\sum_{\pi\in\hat\G} d_u^\pi g(\pi)\int
\overline{\chi_v^\pi(xyx^{-1})}f(xy^{-1})
d\lambda_{u}^u(y)\\
&=\sum_{\pi\in\hat\G} d_u^\pi g(\pi)\int \overline{\chi_u^\pi(y)}f(xy^{-1})d\lambda_{u}^u(y)\\
&=\int \check g(y)f(xy^{-1})d\lambda_{u}^u(y)\\
&=\int \check g(y)f(xy^{-1})d\lambda_{s(x)}(y)\\
&=(f*\check g)(x).
\end{align*}
\qed

\vspace{.3 cm} The diagonal inverse Fourier transform satisfies
the same properties as  the global inverse Fourier transform. We
define natural inner product on $C_c(\hat \G)$ by
$$<g,h>=\sum_{\pi\in\hat\G} {d_u^\pi}^2 \bar g(\pi)h(\pi)\quad (g,h\in C_c(\hat\G)).$$

\begin{prop}
For each $g,h\in C_c(\hat \G)$,

$(i) (\check g)\hat{}=g$,

$(ii) (gh)\check{}=\check g*\check h$,

$(iii) (g^*)\check{}=(\check g)^*$,

$(iv) <\check g,\check h>=<g,h>$.\qed
\end{prop}

Next we characterize the central elements in $L^2(\G)$.

\begin{defn}
Consider $\hat \G$ with the discrete topology.  The spectral
measure $d$ on $\hat\G$ is defined by
$$d(\{\pi\})={d_u^\pi}^2\quad (\pi\in\hat\G).$$
\end{defn}

\begin{defn}
Let
$$\G^{'}=\{x\in\G: s(x)=r(x)\}$$
be the isotropy bundle of $\G$. For each $x\in \G^{'}$, let the
conjugacy class of $x$ be
$$\dot x=\{y^{-1}xy: y\in \G^{s(x)}\},$$
and
$${\G^{'}}^\G=\{\dot x: x\in\G^{'}\}.$$
\end{defn}

Let $q:\G^{'}\to {\G^{'}}^\G$  be the canonical projection,
$x\mapsto \dot x$. It is clear that $\G^{'}$ is a closed (and so
compact) subgroupoid of $\G$. Also ${\G^{'}}^{(0)}=\G^{(0)}=X$.
The Haar system $\{\lambda_u,\lambda^u\}_{u\in X}$ of $\G$
restricted to $\G^{'}$ is a Haar system for $\G^{'}$, which in
turn transfers along $q$ to a Haar system on ${\G^{'}}^\G$, which
we denote by $\{\dot\lambda_u,\dot\lambda^u\}_{u\in X}$.
\begin{defn}

Consider $L^2(\G)$ as a bundle  of Hilbert spaces over $X$, whose
fiber at $u\in X$ is $L^2(\G_u,\lambda_u)$. A function bundle
$f=\{f_u\}_{u\in X}\in L^2(\G)$ is called {\it central} if for
each $u\in X$,
$$f_u(x^{-1}yx)=f_u(y)\quad (x\in \G^u, y\in \G_u^u).$$
We denote the set of all central elements of $L^2(\G)$ by
$\mathfrak CL^2(\G)$.
\end{defn}

\begin{lem}
$L^2( {\G^{'}}^\G)\simeq \mathcal Z(L^2(\G^{'}))=\mathfrak
CL^2(\G^{'})$.
\end{lem}
{\bf Proof} The second  equality could be proved as in Lemma 4.2.
For the first, let $f\in \mathfrak CL^2(\G)$, then $\dot f$
defined by $\dot  f_u(\dot x)=f_u(x)\quad (u\in X, x\in \G^{'})$
is well defined. Also for each $u\in X$,
$$\|\dot f_u\|_2^2=\int | \dot f_u(\dot x)|^2 d\dot\lambda_u(\dot x)
=\int | f_u( x)|^2 d\lambda_u( x)=\|f_u\|_2^2<\infty.$$

Conversely if $g\in L^2({\G^{'}}^\G)$, then $f$ defined by
$f_u(x)=g_u(\dot x)\quad (u\in X, x\in \G^{'})$ is clearly in
$\mathfrak CL^2(\G^{'})$ and fibrewise has the same  $L^2$
norm.\qed

\vspace{.3 cm} All these lead to an alternative version of the
Plancherel theorem, which could be proved similarly by restricting
to central parts.

\begin{thm}{\bf (Palncherel Theorem, diagonal version)}
The diagonal inverse Fourier transform  $\check{}: C_c(\hat\G)\to
\mathfrak CC(\G)$ extends (uniquely) to a bijective linear
isometry $\check{}: L^2(\hat\G)\to L^2({\G^{'}}^\G)$ with inverse
$\hat{}:L^2({\G^{'}}^\G)\to L^2(\hat\G)$.\qed

\end{thm}


\end{document}